\documentclass[11pt]{amsart}
\usepackage{amssymb,latexsym}

\author{Rasul Shafikov, Kaushal Verma}
\title{A local extension theorem for proper holomorphic mappings in ${\bb C}^2$}

\newcommand {\BOX} {\rule{2mm}{2mm}}
\newcommand {\cl} {\mathcal}
\newcommand {\bb} {\mathbb} 
\newcommand {\ra} {\rightarrow}  
\newcommand {\ov} {\overline}
\newcommand {\ti} {\tilde}
\newcommand {\sm} {\setminus}
\newcommand {\ga} {\gamma}
\newcommand {\ep} {\epsilon}
\newcommand {\al} {\alpha}

\numberwithin{equation}{section}
\newtheorem{theorem}{Theorem}[section]
\newtheorem{proposition}{Proposition}[section]
\newtheorem{lemma}{Lemma}[section]

\newtheorem{definition}{Definition}[section]

\date{}

\begin{document}
\keywords{proper mapping, extension, Segre varieties}
\subjclass{Primary 32H40}
\begin{abstract}
Let $f : {\cl D} \ra {\cl D'}$ be a proper holomorphic mapping between
bounded domains ${\cl D}, {\cl D'}$ in ${\bb C}^2$. Let $M, M'$ be
open pieces on ${\partial {\cl D}}, {\partial {\cl D'}}$ respectively
that are smooth, real analytic and of finite type. Suppose that the
cluster set of $M$ under $f$ is contained in $M'$. It is shown that
$f$ extends holomorphically across $M$. This can be viewed as a local
version of the Diederich-Pinchuk extension result for proper mappings
in ${\bb C}^2$.  
\end{abstract} 
\maketitle

\section{Introduction}

\noindent Let ${\cl D}, {\cl D'}$ be smoothly bounded, real analytic
domains in ${\bb C}^2$. With no further assumptions on the domains
such as pseudoconvexity, Diederich and Pinchuk (see \cite{DP1}) show
that any proper holomorphic mapping $f : {\cl D} \ra {\cl D'}$ extends
holomorphically across each point of ${\partial {\cl D}}$. The purpose
of this paper is to propose and prove the following local version of
their result.

\begin{theorem}\label{maintheorem}
Let ${\cl D}, {\cl D'}$ be bounded domains in ${\bb C}^2$ and let 
$f : {\cl D} \ra {\cl D'}$ be a proper holomorphic mapping. 
Suppose that $M, M'$ are open pieces of $\partial {\cl D}, \partial {\cl D'}$ 
respectively such that
\begin{itemize}
\item[(i)] $\partial {\cl D}$ (respectively $\partial {\cl D'}$) is
smooth, real analytic and of finite type (in the sense of D'Angelo) in an open neighborhood of
${\overline M}$ (respectively ${\overline M}')$,
\item[(ii)] the cluster set $cl_f(M) \subset M'$.
\end{itemize}
Then $f$ extends holomorphically across each point on $M$.
\end{theorem}

\noindent Some remarks are in order. First, there is no assumption on
the cluster set of $M'$ under $f^{-1}$. In particular, it is possible
that the cluster set of some $z' \in M'$ contains points near which
$\partial {\cl D}$ may have no regularity at all. Nothing can then be
said about extending $f^{-1}$ across $z'$. The main theorem in
\cite{V} can be considered as a weaker version of Theorem 1.1 since it
was proved with an additional assumption on the cluster set of $M'$ 
under $f^{-1}$. Second, $f$ is not assumed to possess any
apriori regularity, such as continuity, near $M$. Third, a recent result of Diederich-Pinchuk (cf. \cite{DP3}), which is valid for $n \geq 2$, assumes that $f$ is continuous up to $M$ but not apriori proper. Theorem 1.1 on the other hand assumes properness but not the continuity of $f$ on $M$.\\ 

\noindent Let us briefly recall the salient features of the proof in
\cite{DP1}. There are two main steps: first, to extend $f$ as a proper
correspondence and second, to show that extendability as a proper
correspondence implies holomorphic extendability. In step one,
extension as a proper correspondence is shown first for all 
strongly pseudoconvex points on $\partial {\cl D}$. Building on this
they construct proper correspondences, in a step by step approach,
that extend $f$ across weakly pseudoconvex points and even the
exceptional points in the `border' (the set $M\setminus (M^+\cup M^-)$
defined below). Moreover, during all these
constructions, it is essential to do the analogous steps for the
(multivalued inverse) $f^{-1}$. This is the main reason why the local
version of this result cannot be directly derived from \cite{DP1}. As
mentioned above the cluster set of $M'$ may contain points of
$\partial {\cl D}$ with no regularity, and therefore, without
additional assumptions on the cluster set of $M'$ under $f^{-1}$, 
no regularity of $f^{-1}$ near $M'$ can be established. \\

\noindent Many local results on the extendability of holomorphic
mappings and correspondences have been obtained by different authors
and without mentioning the entire list we refer to \cite{B}, \cite{BS}
and \cite{FR} as examples to illustrate the flavor of the techniques
used. One common feature in all these results is a `convexity'
assumption, either geometric or of a function theoretic nature, on $M, M'$. For example, this could be in form of the existence of local plurisubharmonic barriers at points of $M, M'$. So one possible approach to theorem 1.1 would be to show the existence of such barrier functions at points of the one dimensional
strata of the `border'. This seems to be unknown as yet. Thus a
somewhat different approach has to be used to prove
Theorem~\ref{maintheorem}.\\\\ 

\noindent Let $T$ be the set of points on $M$ where its Levi form vanishes. 
Then $T$ can be stratified as $T = T_0 \cup T_1 \cup T_2$, where $T_k$ is a 
locally finite union of smooth real-analytic submanifolds of dimension $k=0,1,2$. 
Denote by $M_s^{\pm}$ the set of strongly pseudoconvex (resp. strongly
pseudoconcave)  
points on $M$. Let $M^{\pm}$ be the relative interior, taken with respect to 
the relative topology on $M$, of ${\ov {M_s^{\pm}}}$. Then $M^{\pm}$ is the set of weakly pseudoconvex (resp. weakly pseudoconcave) points of $M$ and the border $M \sm (M^+ \cup M^-)$ clearly separates $M^+$ and $M^-$. It was shown in 
\cite{DFY} and \cite{DP1} that the stratification $T=T_0 \cup T_1 \cup T_2$ 
can be refined in such a way that the two dimensional strata become maximally 
totally real manifolds. Let us retain the same notation $T_k, k=0,1,2$ for the 
various strata in the refined stratification. Further, let $T_k^+ = M^+ \cap T_k$ for all $k$. Then $T_2^+$ is the maximally totally real strata near which $M$ is weakly pseudoconvex. The set 
\begin{equation}
M_e=\left(M \sm (M^+\cup M^-)\right)\cap (T_1 \cup T_0)
\end{equation}
is the exceptional set. It was shown in \cite{DF1} that 
\begin{equation}
\big( M \sm (M^+ \cup M^-) \big) \cap T_2 \subset M \cap {\hat {\cl D}},
\end{equation}
where ${\hat {\cl D}}$ denotes the holomorphic hull of the domain $\cl D$.
Observe that $M_e \cup T_1^+ \cup T_0^+$ is a locally finite collection of 
real analytic arcs and points. A similar decomposition exists for $M'$. Let 
us note two facts with the hypotheses of Theorem 1.1. First, $f$ clearly
extends across points in $M \cap \hat {\cl D}$. Second, it was shown
in \cite{V} that $f$ also extends across $M^+_s \cup T^+_2$. Thus the
verification of Theorem \ref{maintheorem} will follow once $f$ is
shown to extend across $M_e\cup T^+_0\cup T^+_1$. \\

\noindent Removability of real analytic arcs and points was also considered 
in \cite{DFY} with the assumption that $f : {\cl D} \ra {\cl D'}$ is biholomorphic
and has an extension that is continuous up to $M$. A related result
was obtained in ${\bb C}^n$, $n \geq 2$, (see \cite{H1}) with the
assumption that $M, M'$ are pseudoconvex and $f$ is continuous up to
$M$. \\

\noindent To conclude let us note one consequence of theorem 1.1.

\begin{theorem}
With the hypothesis of theorem 1.1, the extended mapping $f : M \ra M'$ satisfies the additional properties: $f(M^+_s) \subset {M'_s}^+, f(T_2^+) \subset {T'_2}^+, f(T_1^+ \cup T_0^+) \subset {M'_s}^+ \cup {T'_1}^+ \cup {T'_0}^+, f(M_e) \subset M'_e$ and $f(M \cap {\hat {\cl D}}) \subset M' \cap {\hat {\cl D'}}$. Moreover, if $z_0 \in T_0^+$ is an isolated point of $T$, then $f(z_0) \in {T'_0}^+$ is also an isolated point of $T'$.
\end{theorem}

\noindent The conclusion $f(T_1^+ \cup T_0^+) \subset {M'_s}^+ \cup {T'_1}^+ \cup {T'_0}^+$ cannot be strengthened to $f(T_1^+) \subset {T'_1}^+$ as the following example shows. Consider the pseudoconvex domain $\Omega = \{(z, w) : \vert z \vert^2 + \vert w \vert^4 < 1\}$ and the proper mapping $\eta(z, w) : \Omega \ra {\cl B}^2$ from $\Omega$ to ${\cl B}^2$ the unit ball in ${\bf C}^2$ defined as $\eta(z, w) = (z, w^2)$. Clearly $T_1^+ \subset \partial \Omega$ is defined by the real analytic arc $\{(e^{i \theta}, 0)\}$ and such points are mapped by $\eta$ to strongly pseudoconvex points.

\section{Notation and Preliminaries}

\noindent The notion of finite type will be in the sense of D'Angelo which means that none of $M, M'$ can contain positive dimensional germs of complex analytic sets. There are other notions such as finite type in the sense of Bloom-Graham and essential finiteness. The reader is referred to \cite{BER} for definitions and details. However, all these notions are equivalent in ${\bb C}^2$.\\

\noindent Segre varieties have played an important role in the study
of boundary regularity of analytic sets and mappings when the
obstructions are real analytic. The word `analytic' will always mean complex analytic unless stated otherwise. Here are a few of their properties
that will be used in this article. For a more detailed discussion and
complete proofs the reader is referred to \cite{DF2}, \cite{DW} and \cite{BER}. Let us restrict ourselves to ${\bb C}^2$ as the case for $n
> 2$ is no different. We will write $z = (z_1, z_2) \in {\bb C} \times
{\bb C}$ for a point $z \in {\bb C}^2$. \\ 

\noindent Pick $\zeta \in M$ and move it to the origin after a
translation of coordinates. Let $r = r(z, {\ov z})$ be the defining
function of $M$ in a neighborhood of the origin, say $U$, and suppose
that $\partial r/{\partial z_2}(0) \not= 0$. If $U$ is small enough,
the complexification $r(z, \ov w)$ of $r$ is well defined by means of
a convergent power series in $U \times U$. Note that $r(z, \ov w)$ is
holomorphic in $z$ and antiholomorphic in $w$. For any $w \in U$, the
associated Segre variety is defined as 
\begin{equation}
Q_w = \{z \in U : r(z, \ov w) = 0\}.
\end{equation}
By the implicit function theorem $Q_w$ can be written as a graph. In
fact, it is possible to choose neighborhoods $U_1 \Subset U_2$ of the
origin such that for any $w \in U_1$, $Q_w$ is a closed, complex
hypersurface in $U_2$ and  
\begin{equation}
Q_w = \{z = (z_1, z_2) \in U_2 : z_2 = h(z_1, \ov w)\},
\end{equation}
where $h(z_1, \ov w)$ is holomorphic in $z$ and antiholomorphic in
$w$. Such neighborhoods will be called a {\it standard pair} of
neighborhoods and they can be chosen to be polydisks centered at the
origin. Note that $Q_w$ is independent of the choice of $r$. For
$\zeta \in Q_w$, the germ of $Q_w$ at $\zeta$ will be denoted by
${}_{\zeta}Q_w$. Let ${\cl S} := \{Q_w : w \in U_1\}$ be the set of
all Segre varieties, and let $\lambda : w \mapsto Q_w$ be the so-called Segre
map. Then ${\cl S}$ admits the structure of a complex analytic set on
a finite dimensional complex manifold. Consider the complex analytic
set 
\begin{equation}
I_w := {\lambda}^{-1}(\lambda(w)) = \{z : Q_z = Q_w\}. 
\end{equation}
If $w\in M$, then $I_w \subset M$ and the finite type assumption on $M$
forces $I_w$ to be a finite collection of points. Thus $\lambda$ is a
proper map in a 
small neighborhood of each point on $M$. Also, note that $z \in Q_w
\iff w \in Q_z$ and $z \in Q_z \iff z \in M$. We shall also have
occasion to use the notion of the symmetric point that was introduced
in \cite{DP1}. This is defined as follows: for $w$ close enough to
$M$, the complex line $l_w$ containing the real line through $w$ and
orthogonal to $M$ intersects $Q_w$ at a unique point. This is the
symmetric point of $w$ and is denoted by ${}^sw$. It can be checked
that for $w$ outside ${\cl D}$, the symmetric point ${}^sw \in {\cl
  D}$ and vice-versa. Moreover, for $w \in M$, ${}^sw = w$.\\ 

\noindent For $z \in {\partial {\cl D}}$, the cluster set $cl_f(z)$ is
defined as: 
\begin{equation}
cl_f(z) = \{w \in {\bb C}^2 : {\rm there \; exists}\; (z_j)_{j = 1}^{\infty} \subset
{\cl D}, z_j \ra z,\; {\rm such \; that}\; f(z_j) \ra w\} 
\end{equation} 
If $K \subset {\partial {\cl D}}$, then $cl_f(K)$ is defined to be the
union of the cluster sets of all possible $z \in K$.\\ 

\noindent For all the notions and terminology introduced here, we
simply add a prime to consider the corresponding notions in the target
space. For example, ${M'_s}^+$ is the set of all strongly pseudoconvex
points in $M'$ and for $w'$ close to $M'$, $Q'_{w'}$ is the
corresponding Segre variety.\\ 

\noindent Finally, we recall that for an analytic set $A \subset U
\times U' \subset \mathbb C^2\times \mathbb C^2$ of pure dimension two
with  proper projection to the first component $U$, there exists a
system of canonical defining functions 
\begin{equation}\label{canon-eq}
\Phi_I (z,z') = \sum_{|J|\le m} \phi_{IJ}(z) {z'}^J,\ |I|=m, \ 
(z,z') \in U \times U' \subset \mathbb C^2 \times \mathbb C^2.
\end{equation}
Here $\phi_{IJ}(z) \in {\cl O}(U)$ and $A$ is precisely the set of common 
zeros of the functions $\Phi_I(z,z')$. For details see e.g. \cite{C}.
Analytic set $A$ with such properties is usually called a {\it holomorphic
correspondence}. The set  $A$ is called a {\it proper} holomorphic
correspondence if both coordinate projections are proper.


\section{Strategy for the proof of Theorem \ref{maintheorem}}

\noindent As noted earlier, $f$ extends across 
$M \sm (M_e \cup T_1^+ \cup T_0^+)$. 
We will first consider the one dimensional components of $(M_e \cup T_1^+)$. 
So let $\ga$ be a connected, real analytic arc in $M_e \cup T_1^+$ and 
suppose $0 \in \ga$. Let $0 \in U_1 \Subset U_2$ be a standard pair of neighborhoods
small enough so that $f$ extends across $(M \cap U_2) \sm \ga$. This is possible due 
to the fact that $M_e \cup T_1^+ \cup T_0^+$ is a locally finite union of connected 
components. Consider 
\begin{equation}
C := \{w \in U_1 : \ga \cap U_1 \subset Q_w\},
\end{equation}
which is a finite set (see
Lemma 2.3 in \cite{DFY}). Indeed, $\ga \cap U_1 \subset Q_w$ implies
that $Q_w$ is the unique complexification of $\ga \cap U_1$. The Segre
map $\lambda$ is locally proper near the origin and hence the
finiteness follows. 

\noindent We need to to show that $f$ extends: 
(i) across $(\ga \cap U_1) \sm C$ and 
(ii) across the discrete set $\gamma\cap C$. 
The latter case will follow immediately from the former. Indeed, for 
$z \in C$ choose a smooth, real analytic 
arc $\ti \ga \subset M$ containing $z$ such that $\ti \ga$ is transverse 
to $Q_z \cap M$ at $z$. Then the (unique) complexification of $\ti \ga$ 
is distinct from $Q_z$. Therefore the argument of case (i) can be applied 
to prove (ii). Thus without loss of generality we may assume that 
$0 \not \in C$. To show that $f$ extends across the zero dimensional 
strata of $M_e \cup T_1^+ \cup T_0^+$ it suffices to invoke the 
argument of case (ii), which in turn depends on case (i). All the above reductions of the problem can be summarized as follows.

\medskip
\noindent{\bf General Situation:}{\it The set $\gamma$ is a one
  dimensional stratum of $M_e \cup T^{+}_{1}$ and $0\in
  \gamma$. Neighborhoods $U_1$ and $U_2$ are chosen in such a way that
  for all $z \in \ga \cap U_2$, $z \not= 0$, we have $Q_z \not= Q_0$ and $z
  \not \in Q_0$, and therefore, $(Q_0 \setminus \{0\}) \cap M \cap U_2
  \subset M \setminus (M_e\cup T^{+}_{1}) $. Let $U'$ be an open neighborhood in
  ${\bb C}^2$  such that $M'\subset U'$ and $U'$ is small enough so that for $w' \in
U' \sm {\ov {\cl D}}'$, $Q'_{w'}$ is well defined. Does $f$ extend to
  a neighborhood of the origin? 
}
\medskip

\noindent The setup of the General Situation will be henceforth
assumed, unless otherwise stated. Following \cite{DP1} define: 

\begin{definition}  
$V = \{(w, w') \in (U_1 \sm {\ov {\cl D}}) \times (U' \sm {\ov {\cl
      D}}') : f(Q_w \cap {\cl D}) \supset {}_{{}^sw'}Q'_{w'}\}$ 
\end{definition} 

and for all $(w, w') \in V$ consider
\begin{equation}
E = \{z \in {\cl D} \cap U_2 \cap Q_w : f(z) = {{}^sw}', \; f({}_zQ_w)
\supset {}_{{}^sw'}Q'_{w'}\}. 
\end{equation} 
Lemma 12.2 of \cite{DP1} can be applied since $f$ extends across $Q_0
\sm \{0\} \subset M \sm (M_e\cup T^{+}_1)$ to conclude that $V$ is a closed, complex
analytic set in $(U_1 \sm {\ov {\cl D}}) \times (U' \sm {\ov {\cl
    D}}')$ and $\dim \; V \equiv 2$. Furthermore, $V$ has no limit
points on $(U_1 \sm {\ov D}) \times \partial U'$. Following \cite{DP1}
the goal will be to show that $V$ is contained in a closed, complex
analytic set $\hat V \subset (U_1 \sm {\ov {\cl D}}) \times  U'$ with
no limit points on $(U_1 \sm {\ov {\cl D}}) \times \partial U'$. Then
$\hat V$ will define a correspondence that extends $f$ across the
origin and thus we have holomorphic extension as well. What are the
possible obstructions in this strategy? Clearly, the only difficulty is that $V$ may have limit points on $(U_1 \sm {\ov {\cl D}}) \times M'$. So suppose that $(w_0, w'_0) \in {\ov V} \cap ((U_1 \sm {\ov {\cl D}}) \times M')$ and let
$(w_j, w'_j) \in V$ be a sequence of points converging to $(w_0, w'_0)$. Pick
$\zeta_j \in Q_{w_j} \cap U_2 \cap {\cl D}$ such that $f(\zeta_j) =
{}^sw'_j$. Since $w'_j \ra w'_0 \Rightarrow {}^sw'_j \ra {}^sw'_0 =
w'_0$, the properness of $f$ implies that $\zeta_j$ converges to
$\partial {\cl D}$. By Lemma 12.2 of \cite{DP1}, $E \Subset U_2$ and
hence $\zeta_j \ra \zeta_0 \in M \cap U_2$, after perhaps passing to a
subsequence. The following cases arise:\\ 

\noindent (a) If $\zeta_0 \not \in \ga \cap U_2$, then $f$ extends
across $\zeta_0$ and $f(\zeta_0) = w'_0$. Select $z_j \in U_2 \sm {\cl
  D}$ close to $\zeta_0$ such that $f(z_j) = w'_j$ and $z_j \ra
\zeta_0$. The invariance property shows that 
\[
f(Q_{z_j} \cap {\cl D}) \supset {}_{{}^sw'_j}Q'_{w'_j}
\]
for all $j$. Since $(w_j, w'_j) \in V$,
\[
f(Q_{w_j} \cap {\cl D}) \supset {}_{{}^sw'_j}Q'_{w'_j}
\]
for all $j$. Combining these $Q_{w_j} = Q_{z_j}$ and passing to the
limit gives $Q_{w_0} = Q_{\zeta_0}$ and this is a contradiction since
$I_{\zeta_0} \subset M$.\\ 

\noindent (b) If
$w'_0 \in M' \cap {\hat {\cl D}}'$, then Lemma 3.1 in \cite{DP1} shows
that $\zeta_0 \in M \cap \hat {\cl D}$ and therefore $f$ extends
across $\zeta_0$. The same argument as in case (a) yields a contradiction.\\

\noindent (c) If $\zeta_0 \in \ga \cap U_2$, then $f$ is not apriori
known to extend to a neighborhood of $\zeta_0$. The following 
possibilities may occur:
\begin{itemize}
\item[(i)] $cl_f(\zeta_0) \cap {M'_s}^+ \ne \emptyset$,
\item[(ii)] there exists 
$w'_0\in cl_f(\zeta_0) \cap ({T'_2}^+ \cup {T'_1}^+ \cup {T'_0}^+ \cup M'_e)$
\end{itemize}
To prove Theorem \ref{maintheorem} it is enough to consider case
(c). Section 5 contains the proof of the following proposition.

\begin{proposition}\label{pp1}
Let $0 \in (\ga \cap U_1) \sm C$ and suppose that $0' \in cl_f(0) \cap {M'_s}^+$. Then $f$ extends holomorphically to a neighborhood of the origin and $f(0) = 0'$. 
\end{proposition}

\noindent By this proposition, it follows that $f$ will extend across $\zeta_0$ in case (c) (i) and therefore the argument used in case (a) applies again. Hence $V$ cannot have any limit points on 
$(U_1 \sm {\ov {\cl D}}) \times {M'_s}^+$.\\ 

\noindent After that the only remaining possibility is
$\zeta_0\in \gamma \cap U_2$ and 
$w'_0\in cl_f(\zeta_0) \cap ({T'_2}^+ \cup {T'_1}^+ \cup {T'_0}^+ \cup M'_e)$. 
Suppose that $w'_0 \in
{T'_2}^+$. Observe that $\zeta_j \in Q_{w_j}$ for all $j$ and hence
$\zeta_0 \in Q_{w_0}$ in the limit. Therefore $w_0 \in Q_{\zeta_0}$,
where $\zeta_0 \in \ga \cap U_2$. This motivates the consideration of 
\begin{equation}
{\cl L} = \bigcup_{z \in \ga \cap U_1}(Q_z \cap U_1)
\end{equation}
which is a real analytic set in $U_1$ of real dimension 3. ${\cl L}$
is locally foliated by open pieces of Segre varieties at all of its
regular points. Note that ${\cl L} \cap M$ has real dimension at most
two since $M$ is of finite type. Section 6 studies the limit points of
$V$ on ${\cl L} \times {T'_2}^+$ and shows that all of them are
removable. More precisely, we show:

\begin{proposition}\label{pp2}
Let $(p, p') \in {\cl L} \times {T'_2}^+$ be a limit point for $V$. Then ${\ov V}$ is analytic near $(p, p')$.
\end{proposition}

\noindent Finally, as in \cite{DP1}, Bishop's lemma can be applied to
handle the case when $w'_0 \in M'_e \cup {T'_1}^+ \cup {T'_0}^+$. 

\begin{proposition}\label{pp3}
Let $(p, p') \in {\cl L} \times (M'_e \cup {T'_1}^+ \cup {T'_0}^+)$ be a limit point for $V$. Then ${\ov V}$ is analytic near $(p, p')$. As a consequence, $f$ holomorphically extends across the origin.
\end{proposition}

\noindent Thus the proof of Theorem \ref{maintheorem} is complete once
Propositions \ref{pp1} - \ref{pp3} are proved. Theorem 1.2 is
proved in Section 7.\\ 


\section{Properties of the set $V$}

\noindent In this section we prove some technical results concerning
$cl_f(0)$ and derive additional properties of the set $V$. These 
results will be used in the subsequent sections.

\begin{proposition}\label{cluster}
Assume the conditions described in the general situation, and let
$w' \in cl_f(0) \cap M'$. Then there exists a sequence $\{w_j\} \in M
\sm \gamma$, $w_j \ra 0$ such that $f(w_j) \ra w'$. In
particular, we may choose the sequence $\{w_j\} \in M \sm ({\cl L}
\cap M)$. 
\end{proposition}

\noindent {\it Proof:} Fix arbitrarily small neighborhoods $0 \in U$
and $w' \in U'$ and consider the non-empty locally complex analytic
set $A := \Gamma_f \cap (U \times U')$. Without loss of generality, $0
\not \in M \cap {\hat {\cl D}}$. If the proposition were false, the properness of $f$ implies that
$({\ov A} \sm A) \cap (U \times U') \subset (\ga \cap U) \times (M'
\cap U')$. Clearly $(\ga \cap U) \times (M' \cap U')$ is a smooth, real analytic manifold with CR dimension 1. A result of Sibony (cf. theorem 3.1 and corollary 3.2 in \cite{Si}; see \cite{Sk} also) shows that $A$ has locally finite volume in $U \times U'$. Let $\ga^{\bb C}$ be the (local) complexification of
$\ga$. Shrink $U$ if necessary to ensure that $\ga^{\bb C}$ is a
closed analytic set in $U$. Then $A \sm (\ga^{\bb C} \times U')$ is a
non-empty analytic set in $(U \times U') \sm (\ga^{\bb C} \times U')$
with locally finite volume. Now $\ga^{\bb C} \times U'$ is
pluripolar and hence Bishop's theorem (see $\S 18.3$ in cf. \cite{C})
shows that all the limit points of $A \sm (\ga^{\bb C} \times U')$ are
removable singularities. Thus ${\ov A}$ is analytic in $U \times U'$
and ${\ov A} \subset ({\ov {U \cap {\cl D}}}) \times ({\ov {U' \cap
    {\cl D'}}})$.\\

\noindent Next we show that $(0, 0')$ is in the envelope of holomorphy
  of $({\ov {U \cap {\cl D}}}) \times U'$ (cf. \cite{PV}). Let $\pi : {\ov A}
  \ra U$ be the natural projection. Then there are points $z \in
  \pi({\ov A})$ such that $\pi^{-1}(z)$ is discrete. Indeed, if not,
  then each fiber $\pi^{-1}(z)$ is at least 1 dimensional for all $z
  \in \pi({\ov A})$ and hence by theorem 2 in $\S$V.3.2 of \cite{L} it
  follows that 
\[
\dim \; {\ov A} \ge 1 + \dim \; \pi({\ov A}),
\] 
and this implies that $\dim \; \pi({\ov A}) \le 1$. This is a
contradiction since ${\ov A}$ contains $\Gamma_f$ over a non-empty
subset of $U$. For $(z, z') \in {\ov A}$, let $(\pi^{-1}(z))_{(z,
  z')}$ denote the germ of the fiber $\pi^{-1}(z)$ at $(z, z')$. Let 
\[
S = \{(z, z') : \dim \;(\pi^{-1}(z))_{(z, z')} \ge 1\}.
\]
By a theorem of Cartan-Remmert, $S$ is known to be analytic and the
above reasoning shows that $\dim \; S \le 1$. Without loss of
generality $(0, 0') \in S$ as otherwise $\pi : {\ov A} \ra U$ defines
a correspondence near $(0, 0')$ and by theorem 7.4 in \cite{DP1} $f$
would extend across $0$. Select $L \subset {\ov A}$, an analytic
subset in $U \times U'$ with the following properties: it contains
$(0, 0')$, has pure dimension 1 and is distinct from $S$. Note that $L
\subset {\ov {U \cap {\cl D}}} \times U'$. If $L \cap ((M \cap U)
\times U')$ is discrete, then the continuity
principle forces $(0, 0')$ to be in the envelope of holomorphy of
$({\ov {U \cap {\cl D}}}) \times U'$. Since $M$ is assumed to be of
finite type and since $L\not= S$, no open subset of $L$ can be contained in
$L \cap ((M \cap U) \times U')$. The strong disk theorem in \cite{Vl}
shows that $(0, 0')$ is again in the envelope of holomorphy of $({\ov
  {U \cap {\cl D}}}) \times U'$. \\ 

\noindent To conclude, given an arbitrary $g \in {\cl O}(U \cap {\cl
  D})$, we may regard $g \in {\cl O}((U \cap {\cl D}) \times U')$,
i.e., independent of the $z'$ variables. Then $g$ extends across $(0,
0')$ and the uniqueness theorem shows that the extension of $g$, say
${\ti g}$, is also independent of $z'$. Thus $g$ extends across
$0$. In particular, $f$ extends across $0$ and this is a
contradiction.\\ 

\noindent Finally, we can choose $\{w_j\} \in M \sm ({\cl L} \cap
M)$ because ${\cl L} \cap M$ is nowhere dense in $M$. \hfill{\BOX}\\  

\noindent {\it Remarks:} 
\begin{enumerate}
\item First, the cluster set of the origin may be defined
in two possible ways. One way is considered in Section~2, where the
approaching sequence $\{z_j\} \subset {\cl D}$. The other possibility
is to consider $\{z_j\} \subset M \sm (M_e\cup T^{+}_{1})$. This is
well defined since $f$ extends across $M \sm (M_e\cup
T^{+}_{1})$. Clearly, the former set contains the latter in
general. However, in the special case of  the General Situation of Section~3,
Proposition~\ref{cluster} shows that these two notions of the cluster
set coincide.  

\item If $w'$ in Proposition~\ref{cluster} is also known to belong to
${M'_s}^+$, then each point in the sequence $\{w_j\}$ can be chosen
  from $M_s^+$. Indeed, two cases have to be considered. First, if
  $\ga \subset T^+_1$, then for a small enough neighborhood $0 \in U$
  it follows that $M^+_s \cap U$ is dense in $M \cap U$ and hence the
  sequence $\{w_j\}$ can be chosen to belong to $M^+_s$. Second, if
  $\ga \subset M_e$, then every neighborhood $0 \in U$ contains
  points from both $M^+_s$ and $M \cap {\hat {\cl D}}$. It follows
  that for large $j$, $w_j$ must be strongly pseudoconvex since
  $f(w_j)$ is so; $w_j$ cannot belong to any other strata of $M$ as
  otherwise it will be possible to choose a strongly pseudoconcave
  point near the origin that is mapped locally biholomorphically to a
  strongly pseudoconvex point and this violates the invariance of the
  Levi form. 
\end{enumerate}

\bigskip 

\noindent With $V$ as in Definition~3.1 and $0' \in cl_f(0) \cap {M'_s}^+$, 
pick a standard pair of neighborhoods $0 \in U_1 \Subset U_2$ and $0' \in U'
_1 \Subset U'_2$ so that $M' \cap U'_2$ is strongly
pseudoconvex. Choose another neighborhood $U'$, $0' \in U'
\Subset U'_1 \Subset U'_2$ which has the additional property that for all
$w' \in U' \sm {\ov {\cl D}}'$, $Q'_{w'} \cap M' \Subset U'_1$. To see
this, combine the fact that locally $Q'_{0'} \cap {\ov {\cl D}}' = \{0'\}$ with the holomorphic dependence of $Q'_{w'}$ on $w'$. Shrink
$U'_1, U'_2$ so that for all $w' \in M' \cap U'$, $Q'_{w'} \cap M' =
{w'}$. Proposition~\ref{cluster} shows the existence of $(w_0, w'_0) \in ((M \sm
\ga) \cap U_1) \times (M' \cap U')$ such that $w'_0 = f(w_0)$. The
invariance property 
\[
f(Q_w \cap {\cl D}) \supset {}_{{}^sw'}Q'_{w'}
\]
holds for all $(w, w') \in (U_1 \sm {\ov {\cl D}}) \times (U' \sm {\ov
  {\cl D}}')$ close to $(w_0, w'_0)$ and hence $V \cap ((U_1 \sm {\ov
  {\cl D}}) \times U') \not= \emptyset$. Now  $V \cap ((U_1 \sm {\ov
  {\cl D}}) \times U')$ may have several irreducible
components but we retain only those that contain the extension of
$\Gamma_f$ across points in $(M \cap U_1) \times (M' \cap U')$. Take
the union of all such components and call the resulting analytic set
$V_{loc} \subset (U_1 \sm {\ov {\cl D}}) \times (U' \sm {\ov {\cl D}}')$.

\begin{proposition}\label{v1}
$V_{loc}$ consists of those points $(w, w') \in  
(U_1 \sm {\ov {\cl D}}) \times (U' \sm {\ov {\cl D}}')$ that satisfy:
\[
f(Q_w \cap {\cl D}) \supset {}_{{}^sw'}Q'_{w'} \; {\rm and} \; f({}_{{}^sw}Q_w) \subset Q'_{w'} \cap {\cl D'}.
\]
\end{proposition}
\noindent {\it Proof:} We only need to show that for $(w,w')\in V$ the
second inclusion holds. Define 
\[
W = \{(w, w') \in (U_1 \sm {\ov {\cl D}}) \times (U' \sm {\ov {\cl D}}') : f({}_{{}^sw}Q_w) \subset Q'_{w'} \cap {\cl D'}\}.
\]
Pick $(\zeta_0, f(\zeta_0)) \in ((M \sm \ga) \cap U_1) \times (M' \cap
U')$. Then for all $(w, w') \in (U_1 \sm {\ov {\cl D}}) \times (U' \sm
{\ov {\cl D}}')$ close to $(\zeta_0, f(\zeta_0))$ 
\[
f({}_{{}^sw}Q_w) \subset Q'_{w'} \cap {\cl D'}
\]
holds true. Thus $W$ is non-empty. Denote the irreducible component of $V_{loc}$ which contains the extension of $\Gamma_f$ near $(\zeta_0, f(\zeta_0))$ by $\ti V$. It will suffice to show that $\ti V \subset W$.\\

\noindent A similar argument as in \cite{Sh}, Prop. 3.1 shows that $W$ is locally
analytic in $(U_1 \sm {\ov {\cl D}}) \times (U' \sm {\ov {\cl D}}')$
and $\dim \; W \equiv 2$. To show that $W$ is closed, pick a sequence
$(w_j, w'_j) \in W$ converging to $(w_0, w'_0) \in (U_1 \sm {\ov {\cl
    D}}) \times (U' \sm {\ov {\cl D}}')$. Then 
\[
f({}_{{}^sw_j}Q_{w_j}) \subset Q'_{w'_j} \cap {\cl D'}
\]
holds for all $j$. By construction $Q'_{w'_j} \cap M' \Subset U'_1$
and since $f({}^sw_j) \in Q'_{w'_j} \cap {\cl D'}$ we may pass to the
limit to get 
\[
f({}_{{}^sw_0}Q_{w_0}) \subset Q'_{w'_0} \cap {\cl D'}.
\]
Clearly $\ti V \subset W$ locally near $(\zeta_0, f(\zeta_0))$ and
hence it follows that $\ti V \subset W$. \hfill{\BOX}\\ 

\noindent {\it Remarks:} First, without apriori regularity of $f$ near
$0$, such as continuity, it seems difficult to control the images
$f({}_{{}^sw}Q_w)$ as $w$ moves in $U_1 \sm {\ov {\cl D}}$. In
particular, they may escape from $U'_2$. Thus $W$ may not be closed in
general. Second, the two conditions specified in the above proposition
are not the same since $Q_w \cap {\cl D}$ may have many components in
general .

\begin{proposition}\label{v2} 
The set $V_{loc}$ satisfies the following properties:
\begin{enumerate}
\item[(i)]
$V_{loc}$ has no limit points on $(M \cap U_1) \times (U' \sm {\ov {\cl D}}')$,
\item[(ii)]
$V_{loc}$ has no limit points on $(U_1 \sm {\ov {\cl D}}) \times (M' \cap U')$, 
\item[(iii)]
for $w_0 \in (M \sm \ga) \cap U_1$, $\# \pi^{-1}(w_0) \cap ({\ov V_{loc}}
\sm V_{loc}) \cap \left( (M \sm \ga) \times (M' \cap U')\right) \le 1$. 
\end{enumerate}
\end{proposition}

\noindent {\it Proof:} {\it (i)} Suppose that $(w_0, w'_0) \in (M \cap
U_1) \times (U' \sm {\ov {\cl D}}')$ is a limit point for $V_{loc}$. There
are two cases to consider.\\ 

\noindent Case (a) : Let $(w_0, w'_0) \in ((M \sm \ga) \cap U_1)
\times (U' \sm {\ov {\cl D}}')$ and choose a sequence $(w_j, w'_j) \in
V_{loc}$ converging to $(w_0, w'_0)$. Then 
\[
f({}_{{}^sw_j}Q_{w_j}) \subset Q'_{w'_j} \cap {\cl D'}
\]
holds for all $j$. Note that $w_j \ra w_0 \Rightarrow {}^sw_j \ra
{}^sw_0 = w_0$. Since $f({}^sw_j) \subset Q'_{w'_j} \cap {\cl D'}$ and
$Q'_{w'_j} \cap M' \Subset U'_1$, it follows that 
$f({}^sw_j) \ra \zeta'_0 \in Q'_{w'_0} \cap M' \cap U'_1$ 
after perhaps passing to a subsequence. Also, $f$ extends across $w_0$ 
and $f(w_0) = \zeta'_0$. By the invariance property applied near 
$(w_0, \zeta'_0)$ we get 
\[
f({}_{{}^sw_j}Q_{w_j}) \subset Q'_{f(w_j)} \cap {\cl D'}
\]
for all $j$. Combining these inclusions gives $Q'_{f(w_j)} =
Q'_{w'_j}$ for all $j$ and hence $Q'_{\zeta'_0} = Q'_{w'_0}$ in the
limit. This is a contradiction to $I_{\zeta'_0} \subset M'$.\\ 

\noindent Case (b) : Let $(w_0, w'_0) \in (\ga \cap U_1) \times (U'
\sm {\ov {\cl D}}')$ and choose a sequence  $(w_j, w'_j) \in V_{loc}$
converging to $(w_0, w'_0)$. Then 
\[
f(Q_{w_j} \cap {\cl D}) \supset {}_{{}^sw'_j}Q'_{w'_j}
\]
holds for all $j$. Select $\zeta_j \in Q_{w_j} \cap {\cl D}$ such that
$f(\zeta_j) = {}^sw'_j$ and $f({}_{\zeta_j}Q_{w_j}) \supset
{}_{{}^sw'_j}Q'_{w'_j}$. Lemma 12.2 in \cite{DP1} shows that
$\{\zeta_j\} \Subset U_2$ and hence $\zeta_j \ra \zeta_0 \in Q_{w_0}$
after passing to a subsequence. By continuity 
\[
f({}_{\zeta_0}Q_{w_0}) \supset {}_{{}^sw'_0}Q'_{w'_0}.
\]
Note that $\zeta_0 \in {\cl D}$ since ${}^sw'_0 \in {\cl D'}$ and $f$ is proper.\\

\noindent Let us show that $(w_0, w'_0)$ is isolated in $\pi^{-1}(w_0)
\cap ({\ov V_{loc}} \sm V_{loc})$, where $\pi$ is the projection to $U_1$. Suppose
not. Then there exist infinitely many $\{{\ti w}'_j\}$ such that
$(w_0, {\ti w}'_j) \in {\ov V_{loc}} \sm V_{loc}$ and ${\ti w}'_j \ra w'_0$. The
same argument as above can be repeated for each ${\ti w}'_j$, i.e.,
selecting $\zeta_{0, j} \in Q_{w_0}$ etc., so that in the limit we
have
\begin{equation}
f({}_{\zeta_{0, j}}Q_{w_0}) \supset {}_{{}^s{\ti w}'_j}Q'_{{\ti w}'_j}.
\end{equation}
for all $j$. Note that $Q_{w_0} \cap U_2$ will have only finitely many
components, after perhaps shrinking $U_2$ slightly, while (4.1) shows
that at least one component of $Q_{w_0} \cap {\cl D}$ must be mapped
onto infinitely many $Q'_{{\ti w}'_j}$. This contradicts the
properness of the Segre map $\lambda' : U' \ra {\cl S'}$ and the claim
follows.\\ 

\noindent Without loss of generality, we may assume that $w_0 \not \in
M \cap {\hat {\cl D}}$ as otherwise the same argument in case (a)
applies to yield a contradiction. Since $w'_0$ is isolated in the
fiber of ${\ov V_{loc}} \sm V_{loc}$ over $w_0$, it is possible to choose small
balls $B_{\ep}, B'_{\ep}$ centered at $w_0, w'_0$ respectively and of
radius $\ep > 0$ so that the projection 
\begin{equation}
\pi : V_{loc} \cap \big( (B_{\ep} \cap (U \sm {\ov {\cl D}})) \times
B'_{\ep} \big) \ra B_{\ep} \cap (U \sm {\ov {\cl D}}) 
\end{equation}
is proper, and hence a finite branched covering. The canonical
defining pseudo-polynomials of this cover defined as in (\ref{canon-eq})
are monic in $z'_1, z'_2$ 
with coefficients that are holomorphic in $B_{\ep} \cap (U \sm
{\ov {\cl D}})$. Since $w_0 \not \in M \cap {\hat {\cl D}}$,
Trepreau's theorem (cf. \cite{T}) shows that all the coefficients
extend to all of $B_{\ep}$, with a slightly smaller $\ep$
perhaps. Hence there exists an analytic set $V^{ext} \subset B_{\ep}
\times B'_{\ep}$ that contains $V_{loc}$ near $(w_0, w'_0)$. Note that the
projection 
\begin{equation}
\pi : V^{ext} \ra B_{\ep}
\end{equation}
is still proper. Let $(w_{\nu}, w'_{\nu}) \in V_{loc}$ be a sequence that
converges to $(\ti w, {\ti w}') \in V^{ext}$ where ${\ti w} \in (M \sm \ga) \cap U_1$. Such a choice exists by (4.2)
and (4.3). But this is precisely the situation of case (a) and the
same arguments there lead to a contradiction.\\ 

\noindent{\it (ii)} Suppose that $(w_0, w'_0) \in (U_1 \sm {\ov {\cl
    D}}) \times (M' \cap U')$ is a limit point for $V_{loc}$. Select $(w_j,
w'_j) \in V_{loc}$ converging to $(w_0, w'_0)$. Then  
\[
f({}_{{}^sw_j}Q_{w_j}) \subset Q'_{w'_j} \cap {\cl D'}
\]
holds for all $j$. Observe that $\{f({}^sw_j)\} \Subset U'_1 \cap {\cl
  D'}$ since $f$ is proper. On the other hand, $f({}^sw_j) \in
Q'_{w'_j} \cap {\cl D'}$ and as $j \ra \infty$, $Q'_{w'_j} \cap {\cl
  D'} \ra Q'_{w'_0} \cap {\cl D'} = \emptyset$ since $w'_0$ is a
strongly pseudoconvex point. This is a contradiction.\\ 

\noindent {\it (iii)} Finally, suppose that $(w_0, w'_{0, 1}) \in ({\ov V_{loc}} \sm V_{loc}) \cap 
((M \sm \ga) \times (M' \cap U'))$
. Using the same argument as in case (a), it can be seen that $f(w_0)
= w'_{0, 1}$. If $(w_0, w'_{0, 2}) \in ((M \sm \ga) \times (M' \cap
U'))$ is another limit point for $V_{loc}$, then the same argument would
show that $f(w_0) = w'_{0, 2}$ and this is a
contradiction. \hfill{\BOX}\\ 

\noindent {\it Remark:} $V_{loc} \subset (U_1 \sm {\ov {\cl D}}) \times (U'
\sm {\ov {\cl D}}')$ can now be regarded as an analytic set in $U_1
\times (U' \sm {\ov {\cl D}}')$.\\ 


\section{$f$ extends if $cl_f(\zeta_0)\cap {M'_s}^+ \not= \emptyset$} 

\noindent In this section we prove Proposition \ref{pp1}. For that we
will consider the sequence of points $p_j \to 0$, $f(p_j)\to 0'$,
whose existence is guaranteed by Proposition 4.1 and study a certain
family of analytic sets $\{ {\cl C}_{p_j} \}$ associated with
$\{p_j\}$. The goal is to derive some properties of the limit set of
${\cl C}_{p_j}$. We prove several preparation lemmas first.\\

\noindent For any $z \in M \cap U_1$, $Q_z \times (U' \sm {\ov {\cl
    D}}')$ is analytic in $U_1 \times (U' \sm {\ov {\cl D}}')$ of pure
dimension 3. Since $V_{loc}$ contains the extension of the graph of $f$
across some points close to $(0, 0')$, it follows that $V_{loc} \cap (Q_z
\times (U' \sm {\ov {\cl D}}'))$ is either empty or analytic of pure
dimension 1. By Proposition~\ref{cluster}, pick $a \in M \cap U_1$, a strongly
pseudoconvex point across which $f$ extends such that $f(a) \in M'
\cap U'$. Shifting $a$ slightly, if needed, ensures that $a \not \in
     {\cl L} \cap M$. By the invariance property 
\[
f({}_aQ_a) \subset {}_{f(a)}Q'_{f(a)}.
\]
Since both $a, f(a)$ are strongly pseudoconvex, the germs ${{}_a}Q_a$,
${{}_{f(a)}}Q'_{f(a)}$ are contained outside ${\cl D}$ and ${\cl D'}$
respectively. For simplicity we consider representatives of the germs
of ${{}_a}Q_a$ and ${{}_{f(a)}}Q'_{f(a)}$, that satisfy the above
properties. Choose $b \in {{}_a}Q_a \sm {\ov {\cl D}}$
where $f$ is defined, so that $f(b) \in {{}_{f(a)}}Q'_{f(a)} \sm
{\ov {\cl D}}'$. Consider the graph of the extension of $f$ over
${}_{b}Q_a$. This is a pure 1 dimensional germ contained in $V_{loc} \cap
(Q_a \times (U' \sm {\ov {\cl D}}'))$. 
Let 
\begin{equation}
{\cl C_a}\subset V_{loc} \cap (Q_a \times (U' \sm {\ov {\cl D}}'))
\end{equation}
be the irreducible component of dimension 1 that contains this germ. 
Note that ${\cl C_a}$ is
analytic in $U_1 \times (U' \sm {\ov {\cl D}}')$. Also, the invariance
property shows that 
\begin{equation}
{\cl C_a} \subset V_{loc} \cap \big( U_1 \times (Q'_{f(a)} \sm \{f(a)\}) \big).
\end{equation}

\begin{lemma}\label{prep1}
${\ov {\cl C}}_a$ is analytic in $U_1 \times U'$.
\end{lemma}

\noindent {\it Proof:} By Proposition~\ref{v2}, all the limit points of 
${\cl C_a}$ are contained in $(Q_a \cap M \cap U_1) \times \{f(a)\}$.
Since $a \not \in {\cl L} \cap M$, it follows that 
$Q_a \cap M \subset M \sm \ga$ and thus $f$ extends to a neighborhood 
of $Q_a \cap M$. Suppose that 
$(w_0, f(a)) \in {\ov {\cl C}}_a \sm {\cl C_a}$. 
Now exactly the same arguments used in case (a) of Proposition~\ref{v2} 
show that $f(w_0) = f(a)$. The global correspondence 
$f^{-1} : {\cl D'} \ra {\cl D}$ 
has finite multiplicity and hence there can be only finitely many $w_0$ 
so that $(w_0, f(a))$ is a limit point for ${\cl C_a}$. Since 
$\dim \; {\cl C_a} \equiv 1$, the Remmert-Stein theorem shows that 
${\ov {\cl C}}_a$ is analytic in $U_1 \times U'$. \hfill{\BOX}\\

\noindent {\it Remarks:} 
First, it is clear that a different choice of $b \in {{}_a}Q_a$ will 
give rise to the same component ${\cl C_a}$. 
Second, it follows by construction that $(a, f(a)) \in {\ov {\cl C}}_a$ 
and ${\ov {\cl C}}_a \subset Q_a \times Q'_{f(a)}$.\\

\noindent We will now focus on $cl_g(0')$ near $0$, where $g=f^{-1}$ is 
a proper holomorphic correspondence. If $0 \in cl_g(0')$ is not isolated, 
it must be a continuum near $0$. Note that there can only be finitely many $z \in cl_g(0') \cap (M \sm \ga) \cap U_1$ since $g$ has finite multiplicity. Hence we may assume that $cl_g(0') \cap \ga$ is a continuum in $U_1$ containing $0$ and that no point in it belongs to $M \cap {\hat {\cl D}}$. Fix $p \in cl_g(0') \cap (M \sm \ga) \cap U_1$ and a sequence 
$\{p_j\} \subset M^+_s \sm {\cl L}$ converging to $p$ such that $f(p_j)$ converges to $0'$. Associated with each $p_j$ is the analytic set 
${\cl C}_{p_j} \subset U_1 \times (U' \sm {\ov {\cl D}}')$ 
constructed as above such that 
${\ov {\cl C}}_{p_j} \subset U_1 \times U'$ is analytic. The goal will be 
to associate a pure 1 dimensional analytic set, say ${\cl C}_p$, to the chosen point $p$ by considering the sequence ${\cl C}_{p_j}$. 

\begin{lemma}\label{prep2}
Fix $p \in cl_g(0') \cap \ga \cap U_1$ and consider the sequence 
of analytic sets $\{{\cl C}_{p_j}\}$. 
Then the limit set of this sequence of analytic sets is non-empty in $U_1 \times (U' \sm {\ov {\cl D}}')$.
\end{lemma}

\noindent {\it Proof:} Without loss of generality we may assume that
$p=0$. The following observations can be made. First, ${}_0Q_0 \sm
\{0\}$ intersects both $U_1 \cap {\cl D}$ and $U_1 \sm {\ov {\cl D}}$
or else is contained in $U_1 \sm {\ov {\cl D}}$, as otherwise the
continuity principle forces $0 \in M \cap {\hat {\cl D}}$. Thus it is
possible to fix a ball $B_{\ep}$ around the origin so that $(Q_0 \cap
\partial B_{\ep}) \sm {\ov {\cl D}} \not= \emptyset$. Since $Q_z$
depends smoothly on $z$, the same will then be true with $Q_z$ for $z
\in M$ close to $0$. Second, fix a polydisk $\Delta^2 \Subset U'$
centered at $0'$ with its sides parallel to those of $U'$. Then
$Q'_{0'} \cap \partial \Delta^2 \Subset U' \sm {\ov {\cl D}}'$ and the
same will be true for $Q'_{z'}$, $z' \in M'$ close to $0'$. Consider
the non-empty analytic sets ${\ov {\cl C}}_{p_j} \cap (B_{\ep} \times
\Delta^2)$  and examine the projection
\begin{equation}
\pi : {\ov {\cl C}}_{p_j} \cap (B_{\ep} \times \Delta^2) \ra B_{\ep}.
\end{equation} 
Case (a) : If $\pi$ is proper for all $j$, then the image $\pi
\big({\ov {\cl C}}_{p_j} \cap (B_{\ep} \times \Delta^2) \big)$ is
analytic and by the remark after Lemma~\ref{prep1}, it follows that
$\pi \big({\ov {\cl C}}_{p_j} \cap (B_{\ep} \times \Delta^2) \big) =
Q_{p_j} \cap B_{\ep}$. Fix a smaller ball $B_{\ep /2}$ around $0$,
choose $w_0 \in (Q_0 \cap \partial B_{\ep /2}) \sm {\ov {\cl D}}$ and
let $w_j \in (Q_{p_j} \cap \partial B_{\ep /2}) \sm {\ov {\cl D}}$
converge to $w_0$, as $j\to\infty$. Since $\pi$ is proper, it is
possible to choose $(w_j, w'_j) \in   
{\ov {\cl C}}_{p_j} \cap (B_{\ep} \times \Delta^2)$. After passing to
a subsequence, $(w_j, w'_j) \ra (w_0, w'_0)$. Since $w_0 \in U_1 \sm
{\ov {\cl D}}$, Proposition~\ref{v2} shows that $w'_0 \in U' \sm {\ov
  {\cl D}}'$. Thus $(w_0, w'_0) \in (U_1 \sm {\ov {\cl D}}) \times (U'
\sm {\ov {\cl D}}')$ is a limit point for the sequence of analytic
sets $\{{\cl C}_{p_j}\}$.\\ 

\noindent Case (b) : For some subsequence, still indexed by $j$, $\pi$ in (4.5) is not proper. Then it is possible to choose $(w_j, w'_j) \in {\ov {\cl C}}_{p_j} \cap (B_{\ep} \times \Delta^2)$ with $w'_j \in Q'_{f(p_j)} \cap \partial \Delta^2 \Subset U' \sm {\ov {\cl D}}'$. By Proposition~\ref{v2}, $w_j \in U_1 \sm {\ov {\cl D}}$ for all $j$. Passing to a subsequence $w'_j \ra w'_0 \in Q'_{0'} \cap \partial \Delta^2 \Subset U' \sm {\ov {\cl D}}'$ and $w_j \ra w_0 \in {\ov B}_{\ep}$. Thanks to Proposition~\ref{v2} again, $w_0 \in U_1 \sm {\ov {\cl D}}$. Thus $(w_0, w'_0) \in (U_1 \sm {\ov {\cl D}}) \times (U' \sm {\ov {\cl D}}')$ is a limit point for the sequence of analytic sets $\{{\cl C}_{p_j}\}$. \hfill{\BOX}\\

\noindent From (2.1) it follows that $r(z, {\ov p}_j)$ is the defining function for $Q_{p_j}$ in $U_1$. Regard $r(z, {\ov p}_j) \in {\cl O}(U_1 \times U')$, i.e., independent of the $z'$ variables and to emphasize this point, we will write it as $r(z, z', {\ov p}_j)$. This way $r(z, z', {\ov p}_j)$ is the defining function for $Q_{p_j} \times (U' \sm {\ov {\cl D}}')$ in $U_1 \times (U' \sm {\ov {\cl D}}')$ and denote $Z_{p_j} = \{(z, z') : r(z, z', {\ov p}_j) = 0\}$.
Now note that $V_{loc} \cap (Q_{p_j} \times (U' \sm {\ov {\cl D}}'))$ 
is a pure 1 dimensional analytic set. This has two consequences. \\

\noindent First, it follows that (for example, see $\S 16.3$ in
\cite{C}) ${\rm log}\; \vert r(z, z', {\ov p}_j) \vert$ is locally
absolutely integrable on $V_{loc}$. Hence it defines a current,
denoted by ${\rm log} \;\vert r(z, z', {\ov p}_j) \vert \cdot
[V_{loc}]$  
in the following way:
\begin{equation}
\langle  {\rm log}\; \vert r(z, z', {\ov p}_j) \vert \cdot {V_{loc}}, \phi \rangle = \int_{V_{loc}} {\rm log}\; \vert r(z, z', {\ov p}_j) \vert \; \phi,
\end{equation}
where $\phi \in {\cl D}^{(4)}(U_1 \times (U' \sm {\ov {\cl D}}'))$,
the space of all smooth, complex valued differential forms of total
degree 4 with compact support in $U_1 \times (U' \sm {\ov {\cl
    D}}')$. \\ 

\noindent Second, let $\{C_m^j\}_{m \ge 0}$ be the various components
of $V_{loc} \cap (Q_{p_j} \times (U' \sm {\ov {\cl D}}'))$, and let
$k_m^j$ be the corresponding (positive) intersection multiplicities
which are constant along $C_m^j$. Note that for a fixed $j$, ${\cl
  C}_{p_j} = C_m^j$ for some $m$. The wedge product 
\begin{equation}
V \wedge Z_{p_j} = \sum_{m \ge 0} k_m^j C_m^j
\end{equation}
is thus well defined as an intersection multiplicity chain. By the generalized Poincar\'{e}-Lelong formula (cf. \cite{C})
\begin{equation}
 V_{loc} \wedge Z_{p_j} = \frac{1}{2 \pi}\; dd^c \big( {\rm log} \;\vert r(z, z', {\ov p}_j) \vert \cdot [V_{loc}] \big) 
\end{equation}
in the sense of currents for all $j$.

\begin{lemma}\label{prep3}
Fix $p \in cl_g(0') \cap \ga \cap U_1$. Let $\{{\cl C}_{p_j}\}$ 
be the sequence as in Lemma \ref{prep2}. Then there exists a 
subsequence $\{{\cl C}_{p_{j_k}}\}$ converging to an analytic 
set ${\cl C}_{p} \subset U_1 \times (U' \sm {\ov {\cl D}}')$. Moreover ${\ov {\cl C}}_p \subset U_1 \times U'$ is also analytic.
\end{lemma}

\noindent {\it Proof:} By Lemma \ref{prep2} ${\cl C}_{p}$ is not empty.
We now show that the sequence $\{{\cl C}_{p_j}\}$ has locally uniformly bounded volume in $U_1 \times (U' \sm {\ov {\cl D}}')$.

\noindent Fix $K \Subset U_1 \times (U' \sm {\ov {\cl D}}')$. Choose a test function $\psi$ in $U_1 \times (U' \sm {\ov {\cl D}}')$, $0 \le \psi \le 1$ such that $\psi \equiv 1$ on $K$. Let $\omega$ be the standard fundamental form on ${\bb C}^4$. Then
\begin{eqnarray*}
{\rm Vol}\;\big( (V_{loc} \wedge Z_{p_j}) \cap K \big) & = & \int_{(V_{loc} \wedge Z_{p_j}) \cap K} \omega\\
& \le & \int_{(V_{loc} \wedge Z_{p_j})} \psi \omega = \frac{1}{2 \pi} \; \int_{V_{loc}} {\rm log}\; \vert r(z, z', {\ov p}_j) \vert \; dd^c(\psi \omega).
\end{eqnarray*}
Since $r(z, z', {\ov w})$ is antiholomorphic in $w$, it follows that ${\rm log}\; \vert r(z, z', {\ov p}_j) \vert \ra {\rm log}\; \vert r(z, z', {\ov p}) \vert$ in $L^1_{loc}$ and hence
\[
\frac{1}{2 \pi} \; \int_{V_{loc}} {\rm log} \; \vert r(z, z', {\ov p}_j) \vert \; dd^c(\psi \omega) 
\lesssim \frac{1}{2 \pi} \; \int_{V_{loc}} {\rm log} \; \vert r(z, z', {\ov p}) \vert \; dd^c(\psi \omega) := C(K, \psi) < \infty
\]
by the dominated convergence theorem. As noted earlier ${\cl C}_{p_j} = C^j_m$ for some $m$ and hence the assertion follows. \\

\noindent By Bishop's theorem (cf. \cite{C}), $\{{\cl C}_{p_j}\}$ has a subsequence that converges to an analytic set ${\cl C}_p \subset U_1 \times (U' \sm {\ov {\cl D}}')$ locally uniformly. ${\cl C}_p$ has pure dimension 1 and will be reducible in general. By Proposition~\ref{v2}, it is known that all the limit points of ${\cl C}_{p_j} \subset (Q_{p_j} \cap M \cap U_1) \times \{f(p_j)\}$ and hence in the limit ${\ov {\cl C}}_p \sm {\cl C}_p \subset (Q_p \cap M \cap U_1) \times \{0'\}$. But $(Q_p \cap M \cap U_1) \times \{0'\}$ is a smooth, real analytic arc that is also pluripolar. Sibony's result combined with Bishop's theorem applied exactly as before shows that ${\ov {\cl C}}_p \subset U_1 \times U'$ is analytic of pure dimension 1.
\hfill{\BOX}\\

\noindent Note that $(p, 0') \in {\ov {\cl C}}_p$. The association of each $p \in cl_g(0') \cap \ga \cap U_1$ with an analytic set is thus complete. Moreover by (4.4) we have
\begin{equation}
{\ov {\cl C}}_p \subset {\ov {V_{loc} \cap \big( U_1 \times (Q'_{0'} \sm \{0'\}) \big)}}
\end{equation}
for all such $p$.

\begin{lemma}\label{prep4}
Let $p_1, p_2$ be the points in $cl_g(0') \cap \ga \cap U_1$ such that
$p_1\not\in I_{p_2}$. Then $\dim \; ({\ov {\cl
    C}}_{p_1} \cap {\ov {\cl C}}_{p_2}) < 1$. 
\end{lemma}

\noindent {\it Proof:} It may apriori happen that for $p_1 \not = p_2 \in cl_g(0') \cap \ga \cap U_1$, $Q_{p_1} = Q_{p_2}$. But the set of such points is at most countable since $\lambda : U_1 \ra {\cl S}$ is proper. Thus the content of this proposition lies in the assertion that for most points in $cl_g(0') \cap \ga \cap U_1$, the associated analytic sets are also distinct.\\

\noindent Arguing by contradiction assume that 
${\cl C} \subset {\ov {\cl C}}_{p_1} \cap {\ov {\cl C}}_{p_2}$ 
is an irreducible component of dimension 1. Since ${\ov {\cl C}}_{p_i}$ 
lies over $Q_{p_i}$ for $i = 1, 2$ it follows that 
$\pi({\cl C}) \subset Q_{p_1} \cap Q_{p_2}$. 
By hypothesis, $Q_{p_1} \cap Q_{p_2}$ is discrete and the irreducibility 
of ${\cl C}$ implies that $\pi({\cl C})$ is a point, say 
$w \in Q_{p_1} \cap Q_{p_2}$. By (4.9) it follows that 
${\cl C} \subset \{w\} \times Q'_{0'}$ and thus
\begin{equation}
{\cl C} = \{w\} \times Q'_{0'}. 
\end{equation}
There are two cases to consider. First, if $w \in M$, then (4.10) would force $V_{loc}$ to have limit points on $(M \cap U_1) \times (U' \sm {\ov {\cl D}}')$ and this contradicts Proposition~\ref{v2}. Second, if $w \in U_1 \sm {\ov {\cl D}}$, then (4.10) shows that
\[
f(Q_w \cap {\cl D}) \supset {}_{{}^sw'}Q'_{w'}
\]
for all $w' \in Q'_{0'} \sm \{0'\}$. This contradicts the injectivity of $\lambda' : U' \ra {\cl S'}$. \hfill{\BOX}\\

\begin{lemma}\label{prep5}
${\ov {V_{loc} \cap \big( U_1 \times (Q'_{0'} \sm \{0'\}) \big)}}$ 
is analytic in $U_1 \times U'$ of pure dimension one.
\end{lemma}

\noindent {\it Proof:} By (4.9) ${V_{loc} \cap \big( U_1 \times (Q'_{0'} \sm \{0'\}) \big)}$ is a non-empty analytic set of pure dimension one. Proposition~\ref{v2} says that all of its limit points are contained in $(M \cap U_1) \times \{0'\}$. Let $(w, 0') \in ((M \sm \ga) \cap U_1) \times \{0'\}$ be a limit point for $V_{loc}$. The same argument as in case (a) of Proposition~\ref{v2} shows that $f(w) = 0'$. Since $g : {\cl D'} \ra {\cl D}$ has finite multiplicity, it follows that there are only finitely many $w \in M \sm \ga$ such that $(w, 0')$ is a limit point for $V_{loc} \cap \big( U_1 \times (Q'_{0'} \sm \{0'\}) \big)$. Each of these is removable by the Remmert-Stein theorem. The remaining limit points are contained in $(\ga \cap U_1) \times \{0'\}$. Sibony's result combined with Bishop's theorem as before show that they are also removable. \hfill{\BOX}\\

\noindent {\it Proof of Proposition 3.1:} The first step is to show that $0$ is an isolated point in $cl_g(0') \cap \ga \cap U_1$. If not, let $\al > 0$ be the Hausdorff dimension of the continuum ${\cl E} := cl_g(0') \cap \ga \cap U_1$. Each $p \in {\cl E}$ is associated with a pure 1 dimensional analytic set ${\ov {\cl C}}_p \subset U_1 \times U'$. By (4.9) it follows that
\begin{equation}
\bigcup_{p \in {\cl E}} {\ov {\cl C}}_p \subset {\ov {V_{loc} \cap \big( U_1 \times (Q'_{0'} \sm \{0'\}) \big)}}.
\end{equation}
By Lemma~\ref{prep4}, the Hausdorff dimension of the left side in (4.11) is at least $\al + 2$, while the right side has Hausdorff dimension 2 by 
Lemma~\ref{prep5}. Thus $\al = 0$ and that is a contradiction. Hence $0 \in {\cl E}$ is isolated. \\

\noindent By shrinking the neighborhoods $U_1, U_2, U'$ if needed, it
is possible to peel off a local correspondence $g' : U' \cap {\cl D'}
\ra U_1 \cap {\cl D}$ from the global inverse $g : {\cl D'} \ra {\cl
  D}$. Note that $cl_{g'}(0') = 0$ and $cl_{g'}(M' \cap U') \Subset M
\cap U_1$ after shrinking $U'$ even further perhaps. The analytic set
$V_{loc} \subset U_1 \times (U' \sm {\ov {\cl D}})$ defined in
Section~4 contains the graph of the extension of $g'$ across points
near $0'$. Also, $V_{loc}$ does not have limit points on $\partial U_1
\times (U' \sm {\ov {\cl D}}')$ because $cl_{g'}(M' \cap U') \Subset M
\cap U_1$. Thus 
\[
\pi' : V_{loc} \ra U' \sm {\ov {\cl D}}'
\]
is proper. The canonical pseudo-polynomials defining this cover are 
monic with coefficients that are holomorphic in $U' \sm {\ov {\cl D}}'$. 
All of them clearly extend to $U'$, after perhaps shrinking $U'$ and this 
shows that $g'$ extends as a correspondence. Theorem 4.1 in \cite{V} shows 
that $f$ also extends as a correspondence near $(0, 0')$ and this is enough
to conclude that $f$ extends holomorphically across the origin by theorem 7.4 in \cite{DP1}. \hfill{\BOX}\\

\noindent Recall the strategy of the proof of
Theorem~\ref{maintheorem} outlined in Section~3. Let  $\zeta_0$ be as
before and $cl_f(\zeta_0) \cap {M'_s}^+ \not= \emptyset$.  
Then by Proposition~\ref{pp1}, $f$ extends holomorphically to a 
neighborhood of $\zeta_0$ and therefore the argument used in case (a)
can be applied to obtain a contradiction. Thus the global analytic set $V$ as in definition 3.1 has no limit points on $(U_1 \sm {\ov {\cl D}}) \times {M'_s}^+$.


\section{Removability of ${\cl L} \times {T'_2}^+$, ${\cl L} \times
  (M'_e \cup {T'_1}^+ \cup {T'_0}^+)$ and the extendability of $f$} 

\noindent We will now focus on the global analytic set 
$V \subset (U_1 \sm {\ov {\cl D}}) \times (U' \sm {\ov {\cl D}}')$ 
defined in Section~3. Note that the neighborhoods $U_1, U'$ are chosen
as described in the general situation. The results of Section~5 show that
$V$ has no limit points on $(U_1 \sm {\ov {\cl D}}) \times
{M'_s}^+$. The goal of this section is to study the limit points of
$V$ on ${\cl L} \times ( M'_e \cup {T'_2}^+ \cup {T'_1}^+ \cup
{T'_0}^+)$ and to show that they are all removable singularities.
We begin with the limit points on ${\cl L} \times {T'_2}^+$.\\ 

\noindent To start with, note that ${\cl L} = \bigcup_{z \in \ga \cap
  U_1}(Q_z \cap U_1)$ is defined by a single real analytic
equation. Indeed, $\ga$ can be `straightened' by a change of variables
so that in the new coordinate system it becomes the $x_1$-axis. This
may destroy all previous normalizations of the defining function $r(z,
{\ov z})$. But nevertheless, it is clear that 
\[
{\cl L} = \{r(z, x_1) = 0\}.
\]
By the theorems of Cartan-Bruhat (see Prop. 14 and 18 in \cite{N}), it follows that ${\cl L}_{sng}$ is contained in a real analytic set of dimension at most two and is defined in $U_1$ by finitely many real analytic equations. \\

\noindent {\it Proof of Proposition 3.2:} Suppose that $p \in {\cl L}_{reg}$ and choose neighborhoods $p \in U_p, p' \in U'_{p'}$ so that ${\cl L} \times {T'_2}^+$ is smooth, real analytic in $U_p \times U'_{p'}$. Since ${\cl L}$ is locally foliated by open pieces of Segre varieties, it follows that the CR dimension of ${\cl L} \times {T'_2}^+$ is 1 near $(p, p')$. Note that $({\ov V} \sm V) \cap (U_p \times U'_{p'}) \subset {\cl L} \times {T'_2}^+$. As $\dim \; V \equiv 2$, it follows by theorem 20.5 in \cite{C} that $V$ has analytic continuation, say $V^{ext} \subset U_p \times U'_{p'}$ that is a closed analytic set, after shrinking these neighborhoods. Let $\pi : V^{ext} \ra U_p$ and $\pi' : V^{ext} \ra U'_{p'}$ be the projections and 
\[
S = \{(z, z') \in V^{ext} : \dim \;_{(z, z')}({\pi'})^{-1}(z') \ge 1\}
\]
just as in Proposition~\ref{cluster}. The condition 
\[
f(Q_{w} \cap {\cl D}) \supset {}_{{}^sw'}Q'_{w'}
\]
for $(w, w') \in V$ forces 
\[
\pi' : V \cap (U_p \times U'_{p'}) \ra U'_{p'}
\]
to be locally proper. Thus $\pi'(V \cap (U_p \times U'_{p'}))$ contains an open subset of $U'_{p'}$. By the same reasoning used in Proposition~\ref{cluster}, it follows that $\dim \; S \le 1$.\\

\noindent Suppose that $(a, a') \in ({\ov V} \sm V) \sm S \subset {\cl L} \times {T'_2}^+$. It is then possible to choose neighborhoods $U_a, U'_{a'}$ so that
\begin{equation}\label{Vproj}
\pi' : V^{ext} \cap (U_a \times U'_{a'}) \ra U'_{a'}
\end{equation}
is proper. Since $({\ov V} \sm V) \subset {\cl L} \times {T'_2}^+$,
(\ref{Vproj}) shows that $\pi'(V \cap (U_a \times U'_{a'}))$ contains
a one sided neighborhood of $a'$, say $\Omega' \subset U' \sm {\ov
  {\cl D}}'$. Clearly, $\partial \Omega'$ contains points from
${M'_s}^+$ and this contradicts the fact that $V$ has no limit points
on $(U_1 \sm {\ov {\cl D}}) \times {M'_s}^+$. Thus ${\ov V} \sm V
\subset S$. But the three dimensional Hausdorff measure of $S$ is zero
and by Shiffman's theorem (cf. \cite{C}), it follows that ${\ov V}$
itself is analytic in $U_p \times U'_{p'}$. Therefore $V^{ext} = {\ov
  V}$.\\ 

\noindent This argument works if $p \in {\cl L}_{reg}$. As observed above, ${\cl L}_{sng}$ is contained in a real analytic set of dimension at most 2 that is defined by finitely many equations. Thus it is possible to proceed by downward induction to conclude that ${\cl L} \times {T'_2}^+$ is removable. \hfill{\BOX}\\

\noindent As a consequence ${\ov V}$ is analytic in $\big( (U_1 \sm {\ov {\cl D}}) \times U' \big) \sm \big( {\cl L} \times (M'_e \cup {T'_1}^+ \cup {T'_0}^+) \big)$.\\

\noindent {\it Proof of Proposition 3.3:} Consider the global analytic set $V \subset (U_1 \sm {\ov {\cl D}}) \times (U' \sm {\ov {\cl D}}')$ of definition (3.1). We will show that even the bigger set $(U_1 \sm {\ov {\cl D}}) \times (M'_e \cup {T'_1}^+ \cup {T'_0}^+)$ is removable for ${\ov V}$. As observed in $\S 3$,
\[
\pi : {\ov V} \sm ({{\cl L} \times (M'_e \cup {T'_1}^+ \cup {T'_0}^+)}) \ra U_1 \sm {\ov {\cl D}}
\]
is proper. Note that $M'_e \cup {T'_1}^+ \cup {T'_0}^+$ is a locally
finite union of real analytic arcs and points and is thus a locally
pluripolar set. But such sets are also globally pluripolar by
Josefson's theorem. Hence it is possible to choose a plurisubharmonic
function on ${\bb C}^4$, say $\phi$, such that $(U_1 \sm {\ov {\cl D}})
\times (M'_e \cup {T'_1}^+ \cup {T'_0}^+) \subset \{ \phi = - \infty
\}$. Also, choose $a \in M \sm \ga$ across which $f$ is known to
extend. Fix a small ball $B \Subset (U_1 \sm {\ov {\cl D}}) \sm {\cl
  L}$ close to $a$ so that $f$ is well defined in $B$. Then ${\ov V}
\sm ({{\cl L} \times (M'_e \cup {T'_1}^+ \cup {T'_0}^+)})$ has no
limit points on $B \times M'$. Indeed, suppose $(w_0, w'_0) \in B
\times M'$ is such a limit point, and let $(w_j, w'_j) \in {\ov V} \sm
({{\cl L} \times (M'_e \cup {T'_1}^+ \cup {T'_0}^+)})$ be a sequence
converging to $(w_0, w'_0)$. Then 
\[
f(Q_{w_j} \cap {\cl D}) \supset {}_{{}^sw'_j}Q'_{w'_j}
\]
holds for all $j$. Choose $\zeta_j \in Q_{w_j} \cap {\cl D}$ so that $f(\zeta_j) = {}^sw'_j$. After passing to a subsequence, $\zeta_j \ra \zeta_0 \in M \cap U_2$ and ${}^sw'_j \ra {}^sw'_0 = w'_0$. Since $B \cap {\cl L} = \emptyset$ and $\zeta_0 \in Q_{w_0}$, it follows that $\zeta_0 \not \in \ga$. A contradiction can now be obtained exactly as in case (a) of $\S 3$. The non-empty analytic set $\big( {\ov V} \sm ({\cl L} \times (M'_e \cup {T'_1}^+ \cup {T'_0}^+)) \big) \sm \{\phi = -\infty\}$ thus satisfies all the hypotheses of Bishop's lemma (cf. $\S 18.2$ in \cite{C}). Hence ${\ov V}$ is analytic in $(U_1 \sm {\ov {\cl D}}) \times U'$. The projection
\[
\pi : {\ov V} \ra U_1 \sm {\ov {\cl D}}
\]
 is still proper and thus ${\ov V}$ defines a correspondence in $(U_1
 \sm {\ov {\cl D}}) \times U'$. The canonical defining equations of
 this correspondence have coefficients that are holomorphic in $U_1
 \sm {\ov {\cl D}}$. Since $0 \notin {\hat {\cl D}}$, by Trepreau's
 theorem all the coefficients extend to $U_1$, perhaps after shrinking
 $U_1$ . Thus ${\ov V}$ extends as a correspondence to $U_1 \times U'$
 and this provides a multivalued extension of $f$. Theorem 7.4 in
 \cite{DP1} now shows that $f$ extends holomorphically across the
 origin. \hfill{\BOX}\\ 


\section {Proof of Theorem 1.2} 

\noindent By theorem 1.1, $f$ extends across each point of $M$. Let us begin by observing the following:  Let $z_0 \in M^+$ and consider $f(z_0)$. Suppose that $\{p'_j\} \subset M'$ is a sequence of points converging to $f(z_0)$ such that each $p'_j$ is a strongly pseudoconcave point. By the invariance property of Segre varieties (see theorems 4.1 and 5.1 in \cite{DP1}) it is possible to choose small neighborhoods $z_0 \in U$ and $f(z_0) \in U'$ so that the global inverse correspondence $g : {\cl D'} \ra {\cl D}$ extends as a correspondence, say $\ti g$ and
\[
{\ti g} : U' \ra U
\]
is proper. Let $\sigma' \subset U'$ be the branching locus of ${\ti g}$. Fix $p'_{j_0} \in U'$ for some large $j_0$. By shifting $p'_{j_0}$ slightly we may assume that it is still strongly pseudoconcave but $p'_{j_0} \not \in \sigma' \cup f(T \cap U)$, where $T$ is the set of Levi flat points on $M$ . Let $g_1$ be a branch of ${\ti g}$ that is well defined near $p'_{j_0}$ as a locally biholomorphic map. Then $g_1(p'_{j_0})$ is clearly strongly pseudoconvex and this contradicts the invariance of the Levi form. This shows that $f(z_0)$ cannot belong to the border between the pseudoconvex and pseudoconcave points on $M'$ nor can $M'$ be pseudoconcave near it. Thus $f(z_0) \in {M'}^+$.\\

\noindent Case (a) : If $z_0 \in M_s^+$, then combining the observation above with theorem 1.1 in \cite{CPS} shows that $f(z_0) \in {M'_s}^+$. Thus $f(M_s^+) \subset {M'_s}^+$.\\

\noindent Case (b) : Let $z_0 \in T_2^+$. We know that $M'$ must be pseudoconvex near $f(z_0)$. Suppose that $f(z_0) \in {M'_s}^+$. Let $r, r'$ be the defining functions of $M, M'$ near $z_0, f(z_0)$ respectively. By the Hopf lemma, $r' \circ f$ is a defining function for $M$ near $z_0$ and hence, for $z$ close to $z_0$, the Levi determinant $\Lambda$ transforms as
\[
\Lambda_{r' \circ f}(z) = \vert J_f(z) \vert^2 \Lambda_{r'}(f(z)),
\]
where $J_f$ is the Jacobian determinant of $f$. Since $f(z_0) \in {M'_s}^+$, it follows that $\Lambda_{r'}(f(z)) \not = 0$. Thus the zero sets of $\Lambda_{r' \circ f}(z)$ and $J_f$ coincide near $z_0$. But $\Lambda_{r' \circ f}(z)$ vanishes precisely on $T_2^+$ and hence $J_f$ is zero on a maximally totally real manifold. Thus $J_f \equiv 0$ near $z_0$ and this contradicts the properness of $f$.\\

\noindent If $f(z_0) \in {T'_1}^+ \cup {T'_0}^+$, then we may argue in the following manner. Choose small neighborhoods $z_0 \in U, f(z_0) \in U'$ so that
\[
f : U \ra U'
\]
is proper. Then $E := f^{-1}(({T'_1}^+ \cup {T'_0}^+) \cap U') \cap U$ has real dimension 1 and hence there exists ${\ti z} \in (T_2^+ \cap U) \sm E$. Clearly then $f({\ti z}) \in {M'_s}^+$. This is not possible by the discussion above. Thus $f(T_2^+) \subset {T'_2}^+$.\\

\noindent Case (c) : Let $z_0 \in T_1^+ \cup T_0^+$ and suppose that $f(z_0) \in {T'_2}^+$. Let $z_0 \in U, f(z_0) \in U'$ be small neighborhoods as before so that $f : U \ra U'$ is proper. Then $E := f^{-1}({T'_2}^+) \cap U$ has real dimension 2 and hence it is possible to choose $a \in E \sm (T_1^+ \cup T_0^+)$. Note that $a \in M_s^+$. Thus a strongly pseudoconvex point is mapped into ${T'_2}^+$ and this contradicts case (a). Thus $f(T_1^+ \cup T_0^+) \subset {M'_s}^+ \cup {T'_1}^+ \cup {T'_0}^+$.\\

\noindent Case (d) : Let $z_0 \in T_0^+$ be an isolated point of $T$ and suppose that $f(z_0) \in {M'_s}^+$. Then $J_f(z_0) = 0$ as otherwise $f$ would locally biholomorphically map $z_0$, which is a weakly pseudoconvex point to $f(z_0) \in {M'_s}^+$. This contradicts the invariance of the Levi form. Choose small neighborhoods $z_0 \in U, f(z_0) \in U'$ so that
\[
f : U \ra U'
\]
is proper. Call $Z_f = \{z \in U : J_f(z) = 0\}$, the branching locus
of $f$. We claim that $Z_f$ intersects both $U \sm {\ov {\cl D}}$ and
$U \cap {\cl D}$. Indeed, firstly $Z_f \not \subset {\ov {U \cap {\cl
      D}}}$ as otherwise the continuity principle would force $z_0 \in
M \cap {\hat {\cl D}}$ which is not possible. Secondly, an open piece
of $Z_f$ cannot lie in $M$ due to the finite type condition and hence
$Z_f \cap M$ has real dimension at most one. Finally, let us show that
$Z_f \cap {\cl D} \not= \emptyset$. If not, then observe that by the invariance
property of Segre varieties (cf. \cite{DP1}), $f$ maps $U \cap {\cl
  D}$ to $U' \cap {\cl D'}$, $U \sm {\ov {\cl D}}$ to $U' \sm {\ov
  {\cl D'}}$ and $M$ to $M'$. That is, $f$ preserves the `sides' of
$M$. The same is also true for $f^{-1} : U' \ra U$. Now choose some
branch of $f^{-1}$, say $g_1$ that maps a fixed but arbitrary $a' \in U' \cap {\cl D'}$ to $a := g_1(a') \in U \cap {\cl D}$. Since $Z_f$ does not enter ${\cl
  D}$ by assumption, it is possible to analytically continue $g_1$ along all paths in $U' \cap {\cl D'}$ to get a well defined mapping, still denoted by
$g_1$ and $g_1 : U' \cap {\cl D'} \ra U \cap {\cl D}$. The analytic
set $\Gamma_f \subset U \times U'$ extends $g_1$ as a correspondence
and by theorem 7.4 in \cite{DP1}, it follows that $g_1$ extends as a
holomorphic mapping to all of $U'$, after perhaps shrinking $U'$
slightly and $g_1 : U' \ra U$. This shows that $f$ has a well defined
holomorphic inverse and hence it must be a biholomorphic
mapping. Hence $Z_f = \emptyset$ which is not possible as $z_0 \in
Z_f$. Thus $Z_f$ must intersect both sides of $M$ near $z_0$.\\ 

\noindent Choose $a \in (Z_f \cap M) \sm \{z_0\}$ and note that $a \in M_s^+$. Thus $f$ is a proper mapping between strongly pseudoconvex hypersurfaces near $a, f(a)$ that branches at $a$. This is not possible since the Segre maps of both $M, M'$ are injective and this forces $f$ to be locally biholomorphic. Cases (a) and (b) also rule out the possibility that $f(z_0) \in {T'_2}^+ \cup {T'_1}^+$. Thus $f(z_0) \in {T'_0}^+$ is also an isolated point of $T'$. \\

\noindent Case (e) : If $z_0 \in M_e$, then $M'$ cannot be
pseudoconvex near $f(z_0)$ as otherwise there would exist a strongly
pseudoconcave point close to $z_0$ that is mapped to a strongly
pseudoconvex point close to $f(z_0)$. For the same reason $M'$ cannot
be pseudoconcave near $f(z_0)$. Also, $f(z_0)$ cannot belong to the
two dimensional strata of the border since it is known to be in the
envelope of holomorphy (cf. \cite{DF1}). Thus $f(M_e) \subset M'_e$.\\ 

\noindent To conclude, the arguments used in Lemma 3.1 in \cite{DP1} show that $f(M \cap {\hat {\cl D}}) \subset M' \cap {\hat {\cl D}}'$. \hfill{\BOX}

\vskip1cm

\noindent Rasul Shafikov \hfill Kaushal Verma\\
\noindent Dept. of Mathematics \hfill Dept. of Mathematics\\
\noindent SUNY, Stony Brook \hfill University of Michigan\\
\noindent NY 11794 \hfill Ann Arbor, MI 48104

\end{document}